\newif\ifpictures
\picturestrue

\documentclass[11pt]{amsart}
\usepackage[T1]{fontenc}
\usepackage{charter}
\usepackage[expert]{mathdesign}
\usepackage{latex_base}
\usepackage{tikz-cd}    % For making diagrams
\usepackage{mathrsfs}
\usepackage{mathtools}
\usepackage{macros}
\usepackage{afterpage}
\usepackage{lipsum} % For side brackets
\usepackage{xcolor}
\usepackage{amsmath}
\usepackage[dash,dot]{dashundergaps}
\usetikzlibrary{patterns}
\usepackage{tikz,ifthen,fp,calc}
\usetikzlibrary{shapes,arrows,shadows,positioning}
\usetikzlibrary{backgrounds}
\title[]{
A note on polynomial maps having fibers of maximal dimension
}
\author{Boulos El Hilany}
\thanks{For this work, the author was supported by the Institute of Mathematics, Polish Academy of Sciences}

\begin{document}

\maketitle

\begin{abstract} 
For any two integers $k,n$, $2\leq k\leq n$, let $f:(\mathbb{C}^*)^n\rightarrow\mathbb{C}^k$ be a generic polynomial map with given Newton polytopes. It is known that points, whose fiber under $f$ has codimension one, form a finite set $C_1(f)$ in $\mathbb{C}^k$. For maps $f$ above, we show that $C_1(f)$ is empty if $k\geq 3$, we classify all Newton polytopes contributing to $C_1(f)\neq \emptyset$ for $k=2$, and we compute $|C_1(f)|$.
\end{abstract}

 \markleft{}
 \markright{}
\section{Introduction}
Polynomial maps are some of the most classical structures that can be endowed on a complex vector space, and whose appearance is frequent in applications (see e.g.~\cite{DSS09,Las10,VJL10}). Investigating their topology is every so often the same as understanding the properties that the preimages can have. Some of the tools that are employed for such an objective exploit polynomials' combinatorial properties (see e.g.~\cite{ZaNe90,VTO94,Zah96,ABLMH00,Est13}).

Our main result, Theorem~\ref{th:main}, follows this direction by classifying all \emph{Newton polytopes} for which any corresponding polynomial map 

\noindent $f=(f_1,\ldots,f_k):(\C^*)^n\rightarrow\C^k$ has a fiber of codimension one. This represents another step towards understanding the polyhedral properties of maps admitting a prescribed topology (c.f.~\cite{EH19count}). Our restriction is to integers $k,n$ satisfying $2\leq k\leq n$, since otherwise the problem becomes trivial. 

Let us allude to the precise statement by starting with the following notations. A polynomial $P\in\C[x_1,\ldots,x_n]$ is a finite linear combination $\sum a_wx^w$ of monomials, such that $x^w = x_1^{w_1}\cdots x_n^{w_n}$, and $w=(w_1,\ldots,w_n)\in\N^n$. The set of all exponents $w$ appearing in $P$, and satisfying $a_w\neq 0$, is called the \emph{support} of $P$. Its convex hull $\mathcal{N}(P)$, is a subset of the non-negative orthant $\R^n_+$ of $\R^n$ called the \emph{Newton polytope}. Given a collection $\Delta$ of convex integer polytopes $\Delta_1,\ldots,\Delta_k\subset\R^n$, we denote by $\C^\Delta$ the space of all polynomial maps $f$ above such that $\mathcal{N}(f_i) \subset \Delta_i$. This collection is \emph{independent} if for any $I\subset\{1,\ldots,k\}$, we have $\dim\sum_{i\in I} \Delta_i \geq |I|$ (we use here the \emph{Minkowski sum} $A+B = \{a+b~|~a\in A,~b\in B\}$). 

A generic fiber under any map $f\in\C^\Delta$ has codimension $k$ for an independent above collection $\Delta$ (\cite[Theorem 11]{Kho16}). Similarly to how dependency is a combinatorial characteristic that becomes related to the topology of polynomial maps, we ask the following question. 

\begin{question}\label{quest:main}
What property does a choice $\Delta$ of $k$ independent integer convex polytopes in $\R^n_+$ needs to have so that a generic map $f\in\C^{\Delta}$ admits fibers of codimension $m\leq k -1$?
\end{question}

The term \emph{generic} here refers to points located outside a certain algebraic hypersurface in $\C^{\Delta}$ (see e.g.~\cite{Ber75,Est13,Kho16}). We have made $(\C^*)^n$ to be the source space here since the question becomes substantially easier for $\C^n$ (usually $f^{-1}\big(f(\underline{0})\big)$ has high dimension in $\C^n\setminus (\C^*)^n$). % 
We answer Question~\ref{quest:main} for the case $m=1$. Let $C_1(f)$ be the set of those points in $\mathbb{C}^k$ whose fiber under $f$ is of codimension $1$.

\begin{theorem}\label{th:main}
Let $k,n\in\N$, $2\leq k\leq n$, and $\Delta$ be a collection of independent integer polytopes $\Delta_1,\ldots,\Delta_k$ in $\R_+^n$. Then, we have $\Co\neq\emptyset$ for any $f$ in some open subset of $\C^\Delta$ if and only if $k = 2$, $\Delta_1\subset L_1$, and $\Delta_2\subset\conv( L_1\cup L_2)$, where $L_1,L_2\subset\R^n$ are two parallel lines, with $\underline{0}\in L_1$ and such that the convex hull $\conv(L_1\cup L_2)$ of their union does not contain integer points in its relative interior.

Moreover, with the above conditions on $\Delta$, we have $|\Co|\leq |\Delta_2\cap \N^n\cap L_2|-1$ for any $f\in \C^\Delta$, with equality if in addition $f$ is generic.
\end{theorem} The proof of the main result proceeds as follows. The fiber $f^{-1}(\kappa)$ of a point $\kappa\in\Co$ contains an algebraic hypersurface in $(\C^*)^n$. This implies that all $f_i-\kappa_i$ share some polynomial as a factor. We prove that all such points $f-\kappa$ form a parametrized variety $\Red(\Delta)$ in $\C^\Delta$, and we compute its dimension. The obstruction part follows by eliminating all $\Delta$ that are either dependent, or for which $\dim \Red (\Delta)$ is too small to transversally intersect some $k$-dimensional linear subspace in $\C^\Delta$ containing $f-\kappa$. Finally, collections of polytopes appearing in Theorem~\ref{th:main} correspond to polynomials constructed in \S\S~\ref{subsec:constr} using a refined analysis. 

This approach can be extend to answer Question~\ref{quest:main} whenever $2\leq m \leq k-1$, and $3\leq k$ if we impose some additional conditions on the equations for the fibers of $C_m(f)$. This is the subject of a future work.

\begin{example}
For any $a_i,b_j\in\C^*$ below, the map $f:(\C^*)_{r,s,t}^3\to\C^2_{u,v}$
\[(r,s,t)\mapsto (a_0+ a_1rst,~b_0+ b_1r+b_2r^2st+b_3r^3s^2t^2)\] satisfies the combinatorial properties of Theorem~\ref{th:main}. Indeed, the line $L_i$ here is directed by $(1,1,1)$, and passes through the point $(i-1,0,0)$. 

If $b_1a_1^2 + b_2(u-a_0)a_1 + b_3(u-a_0)^2\neq 0$, then $f^{-1}(u,v)$ is either empty or forms the curve \[\left\lbrace (r,s,t)~\left|~rst =\frac{(u-a_0)}{a_1},~r=\frac{a_1^2(v - b_0)}{b_1a_1^2 + b_2(u-a_0)a_1 + b_3(u-a_0)^2}\right.\right\rbrace.\] 
Otherwise, we get $(u-a_0)/a_1$ is a root of $P(z): = b_1 + b_2z + b_3z^2$ and $v=b_0 \Leftrightarrow f^{-1}(u,v)=\left\lbrace (r,s,t)~|~rst =(u-a_0)/a_1\right\rbrace$. The latter set has codimension one in $(\C^*)^3$, and if $f$ is generic in $\C^\Delta\cong\C^6$ (i.e. $a_1b_1b_3(b_2^2 - 4b_1b_3)\neq 0$), then there exists two points $z\in\C^*$ for which $P(z)=0$, and $\cdim f^{-1}(a_0 + a_1z,b_0)=1$. Note that the open subset mentioned in Theorem~\ref{th:main} can be chosen to be $\left\lbrace f\in\C^\Delta~|~a_1b_1b_3\neq 0\right\rbrace$.
\end{example}

\section{Proof of Theorem~\ref{th:main}}\label{sec:bigproof} 
In what follows, we fix once and for all a collection $\Delta$ of independent integer polytopes $\Delta_1,\ldots,\Delta_k$ in $\R_+^n$. Since it will cause no loss of generality, assume in what follows that all members of $\Delta$ contain the origin $\underline{0}$ of $\R^n$ as a vertex. For $i=1,\ldots,k$, let $\C^{\Delta_i}\cong\C^{|\Delta_i\cap \N^n|}$ denote the space of all polynomials $P:\C^n\rightarrow\C$ such that $\mathcal{N}(P) \subset \Delta_i$ (a polynomial is thus identified with its coefficients). We use $\C^\Delta$ to denote the space $\C^{\Delta_1}\oplus \cdots \oplus \C^{\Delta_k}$.

\begin{remark}\label{rem:notation}
It is more common to define the space of polynomials using their supports instead (c.f.~\cite{GKZ94}). Our case is a particular instance of this where the support of $f_i$ is $\Delta_i\cap \N^n$. 
\end{remark}

Let $\mathcal{U}$ be an open subset of $\C^\Delta$, satisfying the assumption of Theorem~\ref{th:main}, and let $f$ be any point in $\mathcal{U}$. Then, the fiber $f^{-1}(\kappa)$ of a point $\kappa$ in $\Co$, contains the hypersurface $S=\left\lbrace x\in(\C^*)^n~|~h(x)=0\right\rbrace$ for some $h\in\C[x_1,\ldots,x_n]$ having more than one term. In other words, there exist relatively-prime polynomials $g_1,\ldots,g_k\in\C[x_1,\ldots,x_n]$ such that 
\begin{equation}\label{eq:pol:reduc}
f_i- \kappa_i = h\cdot g_i,~i=1,\ldots,k.
\end{equation}
\begin{fact*}
The number of points $\kappa$ in $\Co$ is finite. 
\end{fact*}

\begin{proof}
On the one hand, the hypersurface $S$ is a component in the set $\Sin J$ of points $x\in(\C^*)^n$ at which the Jacobian matrix of $f$ is singular. Since $\Sin J$ is an algebraic set, it contains finitely-many components $S_1,\ldots,S_q$, such that $f(S_i)\in\Co$, $i=1,\ldots,q$.

On the other hand, the preimages in $(\C^*)^n$ of any two distinct points $\kappa,\kappa'\in\Co$ are disjoint.

Therefore, a finite number of components $S_1,\ldots,S_q$ produces a finite number of its set of images $\Co$ under $f$.
\end{proof} 

Given $f$, all collections $f-y:=(f_1 - y_1,\ldots, f_k - y_k)$ are identified with the set of points in $\C^\Delta$ such that all $|\Delta_1\cap\N^n| + \cdots + |\Delta_k\cap \N^n|$ coordinates are fixed except the $k$ of them corresponding to parameters $ y_1,\ldots,y_k$. This defines a linear sub-space $H_k(f)$ of dimension $k$ in $\C^\Delta$. 

At the same time, consider the set $\Red(\Delta)$ of \emph{reducible polynomials in $\C^\Delta$}. This is the set, containing the collection $f-\kappa$, and consisting of all collections $\varphi = (\varphi_1,\ldots,\varphi_k)$ above such that for $i=1,\ldots,k$, we have $\mathcal{N}(\varphi_i) = \Delta_i$, and $\varphi_i$ is a product as in~\eqref{eq:pol:reduc}. In particular, the $i$-th member $f_i-\kappa_i$ is written as \[\sum a_{w,i}x^w,\] where $w$ runs through all points in $ \Delta_i\cap\N^n$, and $a_{w,i}$ is some complex polynomial $\Phi_{w,i}$ of positive degree in the coefficients appearing both in $h$, and in $g_i$. Therefore, the set $\Red(\Delta)$ is a union of algebraic varieties in $\C^\Delta$, which we will describe in detail shortly. 

The following statement is a direct consequence of above definitions: \emph{The set $\Co\subset\C^k$ is empty if and only if $H_k(f)\cap\Red(\Delta)\subset\C^\Delta$ is empty}.

In what follows, we pick $f\in\mathcal{U}$ so that we have
\begin{equation}\label{eq:codimensions}
\cdim (H_k(f)\cap \Red (\Delta)) = \cdim H_k(f) + \cdim \Red(\Delta).
\end{equation} This choice is possible since there exists a topological stratification of $\Red\Delta$ with respect to its singularities (e.g. a Whitney stratification) such that one can find a perturbation $f\mapsto\tilde{f}\in\mathcal{U}$ so that $H_k(\tilde{f})$ intersects transversally each stratum. Such a perturbation exists by Sard's theorem since $\mathcal{U}$ is open in $\C^\Delta$, and all strata are smooth.

 If $ \delta_i$ denotes the number $|\Delta_i\cap\N^n|$ for $i=1,\ldots,k$, then the space $\C^\Delta$ has dimension $\sum \delta_i=:\delta$. We thus obtain from~\eqref{eq:codimensions} the following relation on the dimensions
\begin{equation}\label{eq:dimensions}
\dim H_k(f)\cap \Red (\Delta) = \dim \Red(\Delta) + k - \delta.
\end{equation} Next, we will provide conditions on the polytopes $\Delta_1,\ldots,\Delta_k$, and values $k,n$, necessary for $\dim H_k(f)\cap \Red (\Delta)$ to be non-negative. For this task, we show how to parametrize $\Red (\Delta)$.

\subsection{Preliminary analysis} Let $\mathbf{b}_1,\ldots,\mathbf{b}_k,\mathbf{c}$ be the respective coefficients' sets of $g_1,\ldots,g_k,h$ indexed by some common order on the lattice $\N^n$. First, observe that $ \Red(\Delta)$ is $\Ima (\mathbf{\Phi})$, where
\begin{align*}\label{eq:mapPhi}
\mathbf{\Phi}:\bigoplus_{i=1}^k\C^{|\mathbf{b}_i|}\oplus\C^{|\mathbf{c}|} &\longrightarrow \bigoplus_{i=1}^k\C^{|\mathbf{a}_i|}, \\
(\mathbf{b}_1,\ldots,\mathbf{b}_k,\mathbf{c}) & \longmapsto \left( \big( \Phi_{w,1}(\mathbf{b}_1,\mathbf{c})\big)_{w\in \Delta_1\cap\N^n},\ldots,\big( \Phi_{w,k}(\mathbf{b}_k,\mathbf{c})\big)_{w\in \Delta_k\cap\N^n} \right).
\end{align*} Indeed, if $\mathbf{a}_1,\ldots,\mathbf{a}_k$ are the sets of coefficients of $\varphi_1,\ldots,\varphi_k$ respectively enumerated using the same above indexing order on $\N^n$, and if $\varphi\in\Red(\Delta)$, then we have 
\[\mathbf{a}_i = \left( \Phi_{w,i}(\mathbf{b}_i,\mathbf{c})\right)_{w\in\Delta_i\cap\N^n},~i=1,\ldots,k.\] Since $\mathbf{\Phi}$ is a well-defined polynomial map, we get 
\begin{equation}\label{eq:main-ineq1}
\dim\Red(\Delta)=\dim\Ima (\mathbf{\Phi})\leq \sum |\mathbf{b}_i| + |\mathbf{c}|.
\end{equation} On the other hand, Equation~\eqref{eq:pol:reduc} implies 
\begin{equation}\label{eq:Newton}
 \Delta_i= \mathcal{N}(\varphi_i) = \mathcal{N}(h) + \mathcal{N}(g_i).
\end{equation} We prove the following result at the end of this section.

\begin{lemma}\label{lem:difr2}
Let $A,B$, and $\Sigma$ be three integer convex polytopes in $\R^n$ such that $\Sigma$ coincides with the Minkowski sum $A + B$. Let $\alpha$, $\beta$, and $\sigma$ denote the numbers $|A\cap\Z^n|$, $|B\cap\Z^n|$, and $|\Sigma\cap\Z^n|$ respectively. Then, we have $\sigma\geq \alpha + \beta -1$, with equality if and only if $\Sigma$ is contained in a line.
\end{lemma} From this result, the relation in~\eqref{eq:Newton} implies 
\begin{equation}\label{eq:k-inequalities}
|\mathbf{c}| + |\mathbf{b}_i| - 1 \leq |\mathbf{a}_i|,~i=1,\ldots,k.
\end{equation} Recalling that $\delta=\dim \C^{\Delta}$, we obtain from the above inequalities the relation $k(|\mathbf{c}|-1) + \sum |\mathbf{b}_i|  \leq \delta$. Using the second statement of Lemma~\ref{lem:difr2}, one can actually refine the above by writing 
\begin{equation}\label{eq:sharp-ineq}
k(|\mathbf{c}|-1) + \sum |\mathbf{b}_i|  = \delta - r,
\end{equation} where $r$ collects all the differences resulting from possible strict inequalities in~\eqref{eq:k-inequalities}. If the value appearing in~\eqref{eq:dimensions} is non-negative, then equations~\eqref{eq:main-ineq1}, and~\eqref{eq:sharp-ineq} give the inequality 
\begin{equation}\label{eq:ineq1}
k(|\mathbf{c}|-2) \leq |\mathbf{c}| - r.
\end{equation} 

In what follows, we will assume that $\dim \Ima\mathbf{\Phi}\geq 0$, and make an analysis on the different cases that we will obtain.

\subsection{Restrictions} Recall that $\Co\neq \emptyset$, the polynomial $h$ from Equation~\eqref{eq:pol:reduc} has at least two terms (i.e. $|\mathbf{c}|\geq 2$), the target space $\C^k$ is at least two-dimensional, and $\Delta$ is an independent collection of $k$ integer polytopes. The following claims hold true.
\begin{enumerate}[label=(\alph*)]

	\item \label{it:r<3} \emph{We have} $ r\leq 2$. 
	
	Indeed, combining $r\geq 3$ with Equation~\eqref{eq:ineq1} leads to the relation $1\leq -(|\mathbf{c}| - 2)(k - 1)$, which does not hold true in our settings.
	
	\item \label{it:0<r} \emph{We have} $1\leq r$.
	
	Lemma~\ref{lem:difr2} shows that if $r=0$, then each of $\Delta_i$ is contained in a line $L_i$. Moreover, since each $\Delta_i$ has the same summand $\mathcal{N}(h)$ (see~\eqref{eq:Newton}), all those lines coincide. This makes dependent the collection $\Delta_1,\ldots,\Delta_k$.	
	
	\item \label{it:k<4} \emph{We have} $k\leq 3$.
	
	From Item~\ref{it:r<3}, at most two Newton polytopes from the set $\Delta_1,\ldots,\Delta_k$ do \emph{not} form a subset of a line. This readily shows that if $4\leq k$, there exists two distinct $i,j\in\{1,\ldots,k\}$ such that $\dim (\Delta_i +\Delta_j) = 1\leq 2$, making dependent the collection $\Delta_1,\ldots,\Delta_k$.

	\item \label{it:c<=2} \emph{We have} $|\mathbf{c}|= 2$.
	
	This follows from Items~\ref{it:0<r},~\ref{it:k<4}, and Equation~\eqref{eq:ineq1}.

\end{enumerate}

From above four items, we proceed with restrictions as follows. Item~\ref{it:c<=2} yields $h(x) = x^\lambda - t_0x^\mu$ for some $t_0\in\C^*$, and $\lambda,\mu\in\N^n$. Hence, the zero-locus $V(h)$ in $(\C^*)^n$ can be expressed as $\{x\in(\C^*)^n~|~x^u=t_0\}$, where $u=\lambda -\mu \in\Z^n$. Replacing $x^u$ by the parameter $t$, the collection of polynomials $f_1,\ldots ,f_k$ can thus be expressed as the finite linear combinations 
\begin{equation}\label{eq:sys:t0}
P_{\underline{0},i} - \kappa_i + \sum_{\alpha\in \N^n} x^\alpha P_{\alpha,i}(t),~i=1,\ldots,k,
\end{equation} where each $\alpha$ is a vector whose direction is different from that of $u$, and each $P_{\alpha,i}$ is a complex univariate polynomial. Considering that it would not bring less restrictions on $\Delta$, we now take $f$ to be so that no two polynomials $P_{\alpha,i}$ above share a root. In addition, if $x\in V(h)\cap (\C^*)^n$, then for every $i=1,\ldots,k$ we have $P_{\alpha,i}(t_0) = 0$ if $\alpha\neq\underline{0}\in\Z^n$. Indeed, since otherwise there would be an $x\in V(h)$ at which one of the equations~\eqref{eq:sys:t0} is non-zero, contradicting $\kappa\in\Co$. 

This discussion concludes that if $\Co\neq\emptyset$, then there exists at most one $i\in\{1,\ldots,k\}$, and at most one $\alpha\neq\underline{0}$ such that the polynomial $P_{i,\alpha}$ appearing in~\eqref{eq:sys:t0} is not identically zero. Hence the collection $\Delta_1,\ldots,\Delta_k$ is dependent if $k=3$ (same reasoning as Item~\ref{it:k<4}).

\subsection{Constructions}\label{subsec:constr} The above discussion eliminates all collections of Newton polytopes except one family up to index permutations. That is $k=2$, and $f_1,f_2$ are written as \[f_1 = P_{\underline{0},1}\in\C[x^u]\text{,~and}~f_2=P_{\underline{0},2} + x^v P_{v,2}\in \C[x^u] + \C[x^u]x^v,\] where $u,v\in \N^n$ are two primitive integer vectors. Therefore, the set $\Co$ is expressed \[\big\{y\in\C^2~\big|~y_1 = P_{\underline{0},1}(t),~y_2=P_{\underline{0},2}(t),~t\in V(P_{v,2})\big\},\] and thus it contains at most $\deg P_{v,2}$ points. Finally, since the degree of a dense univariate polynomial equals the number of its terms minus one, all the statements of Theorem~\ref{th:main} follow. In particular, the equality in the last statement of Theorem~\ref{th:main} follows from the fact that $P_{v,2}$ has exactly $\deg P_{v,2}$ solutions in $\C^*$ for generic $f\in\mathcal{U}$.

\subsection{Proof of Lemma~\ref{lem:difr2}}

Without loss of generality, we may assume that  one of the vertices of $\Sigma$ contains the origin $\underline{0}$ of $\R^n$. This implies that each of $A$ and $B$ also have $\underline{0}$ as a vertex. Consider a generic hyperplane $H\subset\R^n$ passing through $\underline{0}$, and bordering a half-space $H^+$ that contains $A$. The genericity assumption on $H$ implies that any translation of $H$ contains at most one point of $A$, and thus $\underline{0}$ is the only point of $A$ belonging to $H$. Next, we order the points of $A$ in the following fashion. A point $a$ in $A\cap\Z^n$ takes the index $i\in\{0,1,\ldots,\alpha-1\}$ if and only if the translated hyperplane $H_i:=\{a\}+ H$ has $i$ elements of $A$ in the open half $n$-space $H_i^-$ (the one containing $\underline{0}$). In particular, the vertex of $A$ at the origin is written as $a_0$. 

Let $v\in\R^n$ be some normal vector of $H$, directed towards $A$. The above ordering implies that the scalar product $\langle v, a_{i+1} -a_i\rangle$ is strictly positive for $i=0,1,\ldots,\alpha-1$. Hence, orthogonal projections $\pi_v$ of the points $\{a_0,a_1,\ldots,a_{\alpha - 1}\}\subset A$ to the line $L_v$, directed by $v$, satisfy the ordering $\pi_v(a_i)<\pi_v(a_{i+1})$ with respect to some fixed orientation of $L_v$. 

One can thus re-write $\Sigma$ as $\cup_{i=1}^{\alpha-1}\Sigma_i$, where \[\Sigma_{i}:=\Sigma_{i-1}\cup\big(\{a_{i}\} + B\big)\] for $i=1,\ldots,\alpha-1$ and $\Sigma_0:=B$. This implies that $\Sigma_i\subsetneq \Sigma_{i+1}$, and thus making the integer \[ K_{i}:=|\Sigma_{i}\cap \Z^n| - |\Sigma_{i-1}\cap \Z^n|\] to be always positive. Therefore, the number of elements $\alpha$ is computed as $\beta + \sum_{i=1}^{\alpha - 1} K_i$. 

The equality $\sigma=\alpha + \beta -1$ holds in the sole case where $K_i=1$ for $i=1,\ldots,\alpha-1$. This is only possible whenever all the $a_i$'s belong to a line $L_A$, and the set $B$ is contained in a line $L_B$, with both having the same direction up to scalar. This means that $\Sigma$ is a subset of the Minkowski-summed line $L_A+ L_B$, and we are done. 
\subsection*{Acknowledgements} The author thanks the anonymous referee for their valuable suggestions on the earlier drafts, and for pointing out mistakes therein.

The author is grateful to the Mathematical Institute of the Polish Academy of Sciences in Warsaw for their hospitality and financial support during which this work was accomplished. Support by the Austrian Science Fund (FWF): P33003 is also acknowledged.

\subsection*{Contact}
  Boulos El Hilany,
  Johann Radon Institute for Computational and Applied Mathematics,\\
  Altenberger Stra{\ss}e 69,
  4040 Linz
  Austria; \\
\href{mailto:boulos.hilani@gmail.com}{boulos.hilani@gmail.com}.

\bibliographystyle{alpha}					   % For the style

\bibliography{mainbib}

\end{document}